\documentclass[vecphys]{svmult}    

\usepackage[utf8]{inputenc}
\usepackage{makeidx}         
\usepackage{graphicx}        
\usepackage{multicol}        
\usepackage{multirow}
\usepackage{amsmath}
\usepackage{amssymb}
\usepackage{psfrag}
\usepackage{url}
\usepackage[bottom]{footmisc}

\makeindex             

\newcommand{\bit}{\begin{itemize}}

\newcommand{\eit}{\end{itemize}}

\newcommand{\vk}{{\bf {k}}}

\newcommand{\vE}{{\bf {E}}}

\newcommand{\vu}{{\bf{u}}}
\newcommand{\vU}{{\bf{U}}}
\newcommand{\vW}{{\bf{W}}}
\newcommand{\vf}{{\bf{f}}}
\newcommand{\vF}{{\bf{F}}}

\newcommand{\vx}{{\bf {x}}}
\newcommand{\vn}{{\bf {n}}}

\newcommand{\vy}{{\bf {y}}}
\newcommand{\vz}{{\bf {z}}}
\newcommand{\vv}{{\bf {v}}}

\newcommand{\vb}{{\bf {B}}}

\newcommand{\vnabla}{\nabla}

\newcommand{\vB}{{\bf {B}}}

\newcommand{\deps}{\mathrm{d}\epsilon}

\newcommand{\dx}{\mathrm{d}x}
\newcommand{\dy}{\mathrm{d}y}
\newcommand{\dz}{\mathrm{d}z}
\newcommand{\Dx}{\Delta x}
\newcommand{\Dy}{\Delta y}
\newcommand{\Dz}{\Delta z}
\newcommand{\Dt}{\Delta t}

\newcommand{\beqa}{\begin{eqnarray}}
\newcommand{\eeqa}{\end{eqnarray}}
\newcommand{\beq}{\begin{equation}}
\newcommand{\eeq}{\end{equation}}
\newcommand{\pat}{\partial_t}

\usepackage{algorithm}
\usepackage{algorithmic}

\begin{document}
\title*{Higher-order magnetohydrodynamic numerics}
\titlerunning{Higher-order MHD numerics}
\author{Jean-Mathieu Teissier \and Wolf-Christian M\"uller}
\institute{
J.-M. Teissier \and W.-C. M\"uller
\at Plasma-Astrophysik, ER 3-2, Zentrum f\"ur Astronomie und Astrophysik,
Hardenbergstr.~36a, Technische Universit\"ut Berlin, 10623~Berlin, Germany\\
\url{https://www-astro.physik.tu-berlin.de/node/340}}
%
%
\maketitle
\newcommand{\vfa}{m}
\newcommand{\As}{A}
\newcommand{\Ls}{L}

\newcommand{\Dop}{D}

\abstract{
In this chapter, we aim at presenting the basic techniques necessary to go beyond
the widely accepted paradigm of second-order numerics. We specifically focus on 
finite-volume schemes for hyperbolic conservation laws occuring in
fluid approximations such as the equations of
ideal magnetohydrodynamics or the Euler equations of gas dynamics. For the sake of clarity, a
simple fourth-order ideal magnetohydrodynamic (MHD) solver which allows to
simulate strongly shocked systems serves as an instructive example.
Issues that only or mainly arise in the world of higher-order numerics are given specific focus.
Alternative algorithms as well as
refinements and improvements are dicussed and are
referenced to in the literature.
As an example of application, some results on decaying compressible turbulence are presented.
}

\section{Introduction}
The nonlinear, typically multiscale and multi-physics character of, e.g., the turbulent dynamics of
the solar wind,
the solar corona, or planetary magnetospheres,
requires reliable, efficient and accurate numerical simulations as integral part of theoretical research
in basically all domains of space-plasma physics.
In this context, the description of plasmas as single- or multi-fluid systems has proven to be a powerful
approximation. This is particularly true in settings where the scales of interest
are beyond the realm of kinetic plasma modelling.
A reflection of the level of physical simplification entailed in the fluid approximation is the
reduction of the phase-space particle distribution function to its lowest order velocity moments: mass, momentum,
and energy. Neglecting non-ideal dissipative processes, those quantities are strictly conserved.
This assumed ideality of plasma dynamics is a common approximation given the low levels of collisionality
of many space plasmas. As a consequence, astrophysical numerics typically has to deal with nonlinear conservation laws of
the above-mentioned quantities. Those partial differential equations are of hyperbolic type which means
that the character of their solution is mainly determined by the speeds of signal propagation associated
with the underlying physical processes such as fluid advection or wave propagation. 
If the fluid velocity is of the order of the local sound speed or exceeds it -- a rather common situation
in space -- the solution can lose its regularity and propagating discontinuities, shock fronts,
emerge.

The numerical challenge consists of dealing with the local loss of
differentiability while still guaranteeing the conservation properties
of the fluid equations and thus ensuring the physically correct
evolution of the system. Different strategies of spatial
discretization are available to this end, of which we mention only a
few classical examples: the finite-difference approach representing
the physical fields as a regular grid of point values, the
finite-volume way based on volume-averages of the fluid
observables over arbitrarily shaped grid-cells, and the finite-element
technique which generalizes the finite-volume ansatz to a more complex and
flexible framework for the representation of dependent variables and differential equations on a grid cell.

In this contribution as well as in many numerical investigations during the last decades the
finite-volume technique
is chosen as a reasonable compromise between 
flexibility (compared to finite-differences) and complexity (compared to finite-elements).
It evolves the cell-averages of the conserved quantities via their fluxes on the cell boundaries and,
consequently, is conservative by construction.  

\section{General numerical framework} 
\subsection{System of equations to be solved}
\label{sec:equations}
We are interested in solving the ideal MHD equations assuming
adiabatic thermodynamics.  Since we present here a finite-volume
scheme that employs a constrained transport method to evolve the magnetic
field $\vb$, it is convenient to
express the differential evolution laws as two systems of equations. The hydrodynamical variables, mass density $\rho$, velocity $\vv$, and total energy density $e$ (the sum of internal, kinetic and magnetic energy) 
are governed by the continuity equation, the momentum balance and
the energy equation. They are written here in conservative
form:
\beqa
    \label{eq:hydroA}
    \pat \rho  &=& -\vnabla \cdot (\rho \vv),\\
    \pat (\rho \vv) &=& - \vnabla \cdot \left( \rho \vv \vv^T + (p + \frac{1}{2}|\vb|^2)\mathbf{I} - \vb \vb^T \right),\\
    \label{eq:hydroB} \pat  e  &=& -\vnabla \cdot \left((e + p + \frac{1}{2}|\vb|^2)\vv  -  (\vv\cdot\vb)\vb \right),\\\nonumber
\eeqa
with $p$ the thermal pressure and $\mathbf{I}$ the $3\times 3$ identity matrix so that $\vnabla \cdot (p + \frac{1}{2}|\vb|^2)I=\vnabla (p + \frac{1}{2}|\vb|^2)$. Using the adiabatic equation of state for an ideal gas, the internal energy is $\frac{p}{\gamma - 1}$ with $\gamma$ the ratio of specific heats, which means that $p$ can be deduced through:
\beq \label{eq:pressure}
         p =  (\gamma -1) \left(e - \frac{1}{2}\rho |\vv|^2 - \frac{1}{2}|\vb|^2 \right).
\eeq
Please note that in the isothermal case $\gamma=1$ and the pressure is determined by $p=\rho c_s^2$ with 
constant sound speed $c_s$. In that case, no equation for $e$ needs to be solved.

As for the magnetic field, its evolution is governed by:
\beq
	\label{eq:mag}
	\pat \vb = -\vnabla \times \vE,
\eeq
with $\vE = - \vv \times \vb$ the electric field. In addition to the above equations,
the absence of magnetic monopoles requires the solenoidality of the magnetic field at all times:
         \beq \label{divb}
                \vnabla \cdot \vb = 0.
         \eeq
\subsection{Motivation for using higher-order schemes}
\label{sec:motivation}
Numerical schemes whose discretization error is of second-order in the
grid scale are often combined
with more precise Riemann solvers that are required in finite-volume schemes to calculate
the fluxes over the faces of a grid cell (see section \ref{sec:RP}). An alternative 
way to reduce numerical dissipation are higher-order
schemes which increase the order of the overall discretization error to accelerate the convergence towards the
reference solution in continuous space and time as the numerical resolution is increased.
In this chapter, we qualify a scheme as \lq\lq higher-order\rq\rq\ if its discretization order is strictly greater than 3, and of \lq\lq lower-order\rq\rq\ otherwise.
Higher-order schemes by definition represent the solution by considering more terms of the associated
Taylor series at a specific location on the numerical grid. Consequently, they achieve a more accurate
and therefore less dissipative approximation of the reference solution. Although they generally have a
higher computational cost at a given resolution as compared to a
lower-order scheme, a higher-order algorithm allows to obtain the same solution quality at lower numerical
resolution. Therefore and in particular in more than one spatial dimension, their employment allows
a significant reduction of computational cost. Section \ref{sec:numtests} corroborates this statement by comparison of
fourth-order and lower-order schemes.
\subsection{Finite-Volume}
\label{sec:FV}

\begin{figure}
\begin{center}
\includegraphics[width=0.5\textwidth]{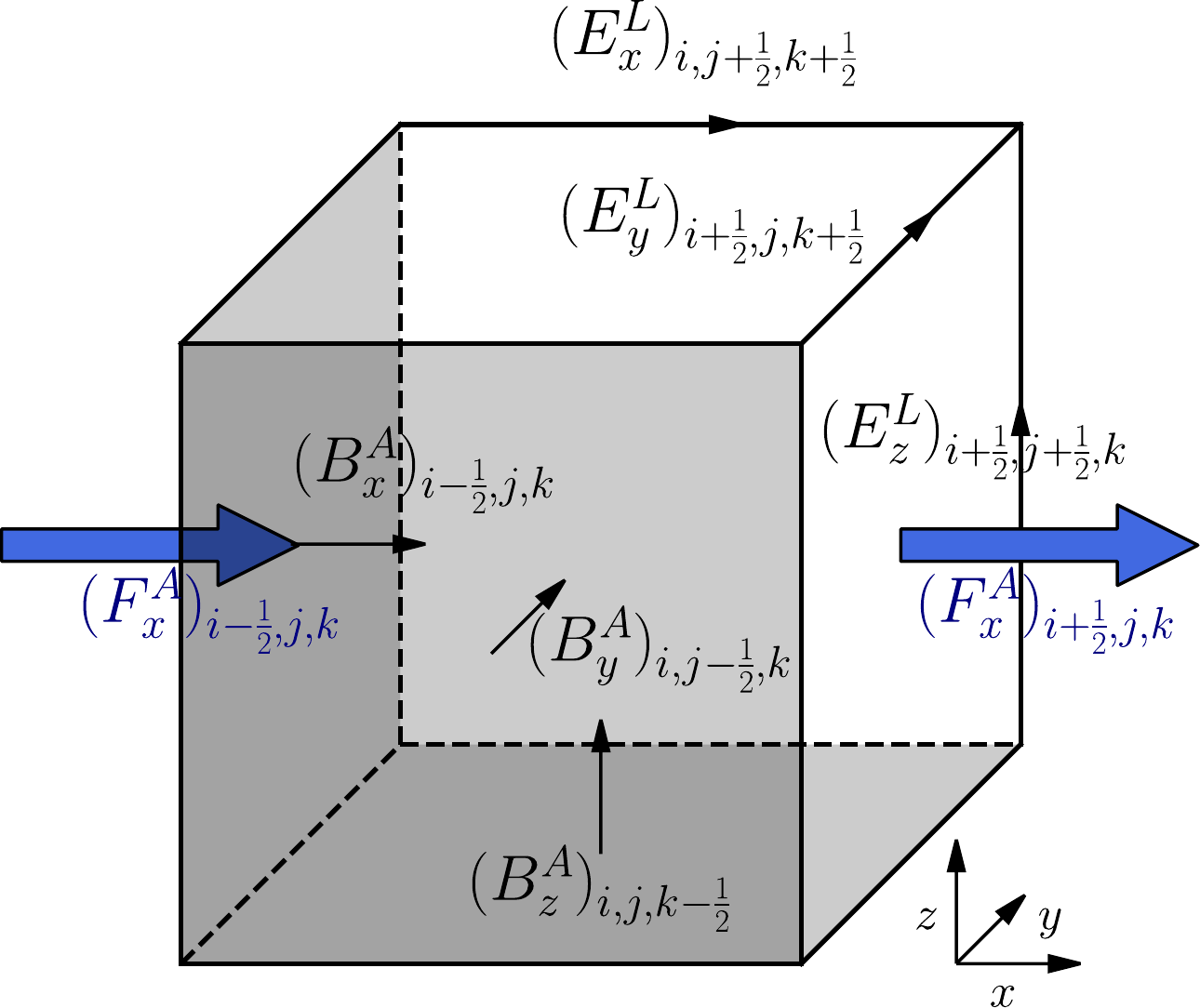}
\caption{ \label{fig:vardef} Definition of the magnetic quantities and the fluxes in the cell $\Omega_{i,j,k}$. Only the flux in the $\vx$-direction is shown, the ones in the other directions are defined analogously.}
\end{center}
\end{figure}

In a finite-volume approach, spatial cell averages are evolved in time. We consider hence the time evolution of the cell average:
\beq
\vU_{i,j,k}=\frac{1}{\Dx\Dy\Dz}\int_{\Omega_{i,j,k}} \vu(x,y,z) {\dx  \, \dy \, \dz \,} ,
\eeq
with $\vu=(\rho,\rho v_x,\rho v_y,\rho v_z,e)$ a vector containing the hydrodynamical quantities and $\vU$ the corresponding cell averages over the cell $\Omega_{i,j,k}=[x_i-\frac{\Dx}{2},x_i+\frac{\Dx}{2}]\times[y_j-\frac{\Dy}{2},y_j+\frac{\Dy}{2}]\times[z_k-\frac{\Dz}{2},z_k+\frac{\Dz}{2}]$. We assume for the sake of simplicity a Carthesian coordinate system where $(x_i=(i+{1/2})\Dx, y_j=(j+{1/2})\Dy, z_k=(k+{1/2})\Dz)$ is the center of the cell indexed by $(i,j,k)$ and $\Dx,\Dy,\Dz$ are constant, not necessarily equal grid-sizes. Using relations \eqref{eq:hydroA}-\eqref{eq:hydroB} and after application of Gauss' theorem, we obtain:
\beq
\label{eq:patvU}
\pat \vU_{i,j,k}=-\frac{\Dop_x(\vF^{{\As},x}_{i,j,k})}{\Dx}-\frac{\Dop_y(\vF^{{\As},y}_{i,j,k})}{\Dy}-\frac{\Dop_z(\vF^{{\As},z}_{i,j,k})}{\Dz},
\eeq
where we introduced in order to have concise notations a difference operator:
\beq
\Dop_x(F_{i,j,k})=F_{i+1/2,j,k}-F_{i-1/2,j,k}
\eeq
for any quantity $F$, and similarly $\Dop_y(F_{i,j,k})=F_{i,j+1/2,k}-F_{i,j-1/2,k}$ and $\Dop_z(F_{i,j,k})=F_{i,j,k+1/2}-F_{i,j,k-1/2}$. The fluxes $\vF^{{\As},x}, \vF^{{\As},y}$ and $\vF^{{\As},z}$ along the $\vx-$, $\vy-$ and $\vz-$directions respectively are area-averaged (see Fig. \ref{fig:vardef}), that is:
\beq
\vF^{{\As},x}_{i\pm 1/2,j,k}=\frac{1}{\Dy\Dz} \int_{{\As}^x_{j,k}} \vf^x(x_{i\pm 1/2},y,z) {\dy \, \dz \,},
\eeq
where ${\As}^x_{j,k}=[y_j-\frac{\Dy}{2},y_j+\frac{\Dy}{2}] \times [z_k-\frac{\Dz}{2},z_k+\frac{\Dz}{2}]$ and similarly for $\vF^{{\As},y}$ and $\vF^{{\As},z}$. The fluxes $\vf^x$, $\vf^y$ and $\vf^z$ can be derived from \eqref{eq:hydroA}-\eqref{eq:hydroB}, for example:
\beqa
\label{eq:fluxfx}
{\vf^x} &=&\begin{pmatrix} \rho v_x \\ \rho v_x^2+p+|\vb|^2/2-B_x^2 \\ \rho v_xv_y-B_xB_y \\ \rho v_xv_z-B_xB_z \\ (e+p+|\vb|^2/2)v_x-B_x(\vv \cdot \vb) \end{pmatrix}.
\eeqa
Please note that even though the right-hand side terms of \eqref{eq:patvU} look very similar to second-order approximations of derivatives, this relation is in fact exact.

This formulation has the advantage of conserving the left-hand-side quantities of \eqref{eq:hydroA}-\eqref{eq:hydroB}
up to machine precision, since what exits one cell enters a neighbouring one and vice-versa. Mass density, momentum and total energy are inherently conserved, which means that all numerical dissipation results in the rise
of internal energy, i.e. a heating of the fluid in the adiabatic case.
\subsection{Constrained Transport}
\label{sec:CT}

In MHD simulations, it is very important to maintain the solenoidality of the magnetic field ($\vnabla \cdot \vB=0$), otherwise unphysical effects such as enhanced transport in the direction orthogonal to the local magnetic field \cite{BRA85} and loss of momentum and energy conservation \cite{BRB80} would occur. In fact, if the magnetic field is evolved in time as a cell average, like the hydrodynamical variables, $|\vnabla \cdot \vB|$ typically grows in time.

Several ways exist which address this issue.
Three well-known techniques are: the Helmholtz-Hodge projection \cite{BRB80} which exactly eliminates
any dilatational component of the magnetic field but does not ensure that the projected field is physically
consistent with the dynamical state of the plasma. This leads to a non-negligible generation
of divergence-free but spurious small-scale magnetic fluctuations \cite{BAK04}.
In contrast, the two other popular approaches namely the 8-wave multiplier method \cite{POW94} and the
Generalized Lagrange Multiplier technique \cite{DKK02} do not annihilate magnetic field divergences but
let the plasma flow advect the erroneous field components preventing their spatial accumulation.
In this context, numerical diffusion is seen as a beneficial effect as it dissipates the advected small-scale
errors in the magnetic field over time.
Here, we focus on the Constrained Transport (hereafter CT) technique \cite{EVH88}. It preserves
magnetic field solenoidality up to machine precision in a 
way analogous to the cell-average conservation in finite-volume schemes.
Refinements of this technique have been developed \cite{BAL04,BAL09}. 

In the CT approach, the magnetic field components, $B_x$, $B_y$ and $B_z$, are not defined as volume-averages like the hydrodynamical variables, but as area averages over the faces of a grid cell that are orthogonal to the $\vx-$, $\vy-$ and $\vz-$
direction respectively (see Fig. \ref{fig:vardef}), for example:
\beq
(B^{\As}_x)_{i \pm 1/2,j,k}=\frac{1}{\Dy\Dz}\int_{{\As}^x_{j,k}} B_x(x_{i \pm 1/2},y,z) {\dy \, \dz},
\eeq
and similarly for $B^{\As}_y$ and $B^{\As}_z$. The face area ${\As}^x_{j,k}$ has been defined in section \ref{sec:FV}. Applying Stoke's theorem to \eqref{eq:mag} leads to the exact relations:
\beqa
\label{eq:dbSA}
\pat (B^{\As}_x)_{i \pm 1/2,j,k}&=&-\frac{\Dop_y((E^{\Ls}_z)_{i \pm 1/2,j,k})}{\Dy}+\frac{\Dop_z((E^{\Ls}_y)_{i \pm 1/2,j,k})}{\Dz}, \\
\pat (B^{\As}_y)_{i,j \pm 1/2,k}&=&\frac{\Dop_x((E^{\Ls}_z)_{i,j \pm 1/2,k})}{\Dx}-\frac{\Dop_z((E^{\Ls}_x)_{i,j\pm 1/2,k})}{\Dz}, \\
\label{eq:dbSB}
\pat (B^{\As}_z)_{i,j,k \pm 1/2}&=&-\frac{\Dop_x((E^{\Ls}_y)_{i,j,k\pm 1/2})}{\Dx}+\frac{\Dop_y((E^{\Ls}_x)_{i,j,k \pm 1/2})}{\Dy},\\\nonumber
\eeqa

where $E^{\Ls}$ is the line-averaged electric field (with $\vE=-\vv \times \vb$), as shown in Fig. \ref{fig:vardef}. For example:
\beq
(E^{\Ls}_z)_{i \pm 1/2,j+1/2,k}=\int_{z=z_k-\Dz/2}^{z_k+\Dz/2} E_z(x_{i\pm 1/2},y_{j+1/2},z) \dz\,.
\eeq
The second-order approximation of $\vnabla \cdot \vB$ is then:
\beq
\vnabla \cdot \vB \approx \frac{\Dop_x((B^{\As}_x)_{i,j,k})}{\Dx}+\frac{\Dop_y((B^{\As}_y)_{i,j,k})}{\Dy}+\frac{\Dop_z((B^{\As}_z)_{i,j,k})}{\Dz}
\eeq
Using \eqref{eq:dbSA}-\eqref{eq:dbSB}, it is straightforward to show that this approximation of $\pat(\vnabla \cdot \vB)$ is conserved, cf. \cite{ZIE04}. Thus, if the magnetic field is initially solenoidal, then its solenoidality is preserved up to machine precision as the field is evolved in time.
\section{Practical Computation of the Fluxes}
The computation of the numerical fluxes over the faces of a grid cell is a fundamental operation
in the finite-volume framework. To this end, a reconstruction step has to be performed, based on the
assumption 
that within a grid cell all physical variables have a polynomial representation up to a given order
of the cell extension, e.g. $\Delta x$. The reconstructed polynomial has to explicitly accomodate to
discontinuities occuring at the grid cell boundaries. Especially for higher-oder reconstruction, where
multiple neighbouring cells have to be taken into account, this becomes important and particularly challenging
if the solution itself can feature physical discontinuities like, e.g., propagating shock-fronts in
high-Mach-number flows.
\subsection{Central Weighted Essentially Non Oscillatory reconstruction}
\label{sec:WENO}
\begin{figure}
\begin{center}
\includegraphics[width=0.5\textwidth]{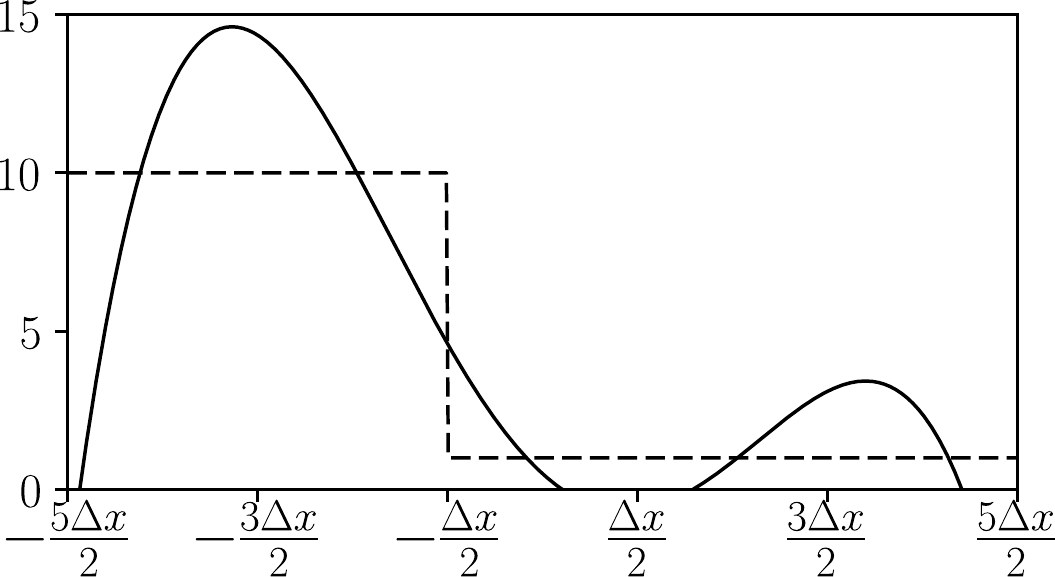}
\caption{ \label{fig:whyweno} Simple polynomial reconstruction over a discontinuity. The dashed line is the function to be reconstructed and the solid line is the unique polynomial of degree 4 which cell averages verify \eqref{eq:poptrel}. The grid-size $\Dx$ is unity.}
\end{center}
\end{figure}

There exist a lot of possibilities in order to compute the area-averaged values $\vU^{\As}_{i \pm 1/2,j,k}$ from the volume-averaged values. For a $(2n+1)^{th}$ order procedure, the most straightforward one is to find for each quantity $q\in\vu$ the unique polynomial $P^{opt,q}_{i,j,k}$ of degree at most $2n$ which cell averages on each of the cells of the $(2n+1)$-cells long stencil $\{x_{i-n},x_{i-n+1},...,x_i,x_{i+1},...,x_n\}$ coincides with $q$'s cell average $Q$:
\beq
\label{eq:poptrel}
\frac{1}{\Dx}\int_{x_{i+\vfa-1/2}}^{x_{i+\vfa+1/2}} P^{opt,q}_{i,j,k}(x) \dx = \int_{\Omega_{i+\vfa,j,k}}q{\dx\,\dy\,\dz}=Q_{i+m,j,k}, \forall \vfa \in [ -n, n ]
\eeq
Then, the $\left(P^{opt,q}_{i,j,k}(x_{i \pm 1/2})\right)_{q\in\vu}$ build a $(2n+1)^{th}$-order approximation of $\vU^{\As}_{i \pm 1/2,j,k}$.
However, using such a polynomial over a discontinuity would be
problematic. As can be seen in Fig. \ref{fig:whyweno}, such a
reconstruction is oscillatory and gives large over- and
undershoots. If the reconstructed value is the density
$\rho$ for example, this would even give rise to an unphysical state with a
negative density at $x=\frac{\Dx}{2}$.

In order to tackle this problem, some
methods use slope limiters which restrict the local variation of the solution to finite
values by introducting numerical dissipation. Examples of
such methods include the Total Variation Diminishing (TVD) limiter of
Van-Leer \cite{VLE77}, a second-order procedure used for example in
the Nirvana code \cite{ZIE04}, the third-order Piecewise Parabolic
Method \cite{COW84} used for example in the
Athena code \cite{SGA08}, the fifth-order MP5 scheme
\cite{SUH97}, etc. Another class of reconstruction methods comprises the
Weighted Essentially Non Oscillatory (WENO) schemes, which we will
present more in detail here. These schemes are improvements of the Essentially Non Oscillatory (ENO) schemes \cite{HEO87}. The first WENO scheme, which was third-order accurate, was designed in 1994 \cite{LOC94}. Since then, very high-order schemes, at least up to order 17, have been derived \cite{BAS00,GSV09}. We present here a refinement of the standard fifth-order WENO \cite{JIS96} used for example in the Pluto code \cite{MBM07}, namely a fourth-order Central Weighted Non Oscillatory (CWENO) scheme \cite{LPR99}. For more literature about WENO schemes, the interested reader can see \cite{SHU97,SHU09}. A comparison between several solvers can be found in \cite{KNC11}.

The idea of this CWENO scheme is the following: instead of using all the cells of a stencil at once to find $P^{opt}$, one divides this stencil into $n+1$ sub-stencils. For example, for $n=2$, the 5-cells wide stencil $\{x_{i-2},x_{i-1},x_{i},x_{i+1},x_{i+2}\}$ is divided in three 3-cells wide stencils containing the central cell: $L=\{x_{i-2},x_{i-1},x_i\}$ on the left, $C=\{x_{i-1},x_{i},x_{i+1}\}$ in the center and $R=\{x_{i},x_{i+1},x_{i+2}\}$ on the right. Then, let $P_L, P_C$ and $P_R$ be the unique polynomials of degree at most 2 which cell averages coincides with the cell averages of $q$ on the stencils $L, C$ and $R$ respectively\footnote{For the sake of simplicity, we drop the superscript $q$ in the following, even though the polynomials are of course different for each variable.}. If $q$ is smooth over these stencils, then each polynomial $P_L$, $P_C$ and $P_R$ gives a third-order approximation of $q$ on the domain $[x_{i-1/2},x_{i+1/2}]$. One can find three positive coefficients $c_L, c_C$ and $c_R$ with $c_L+c_C+c_R=1$ such that:

\beq
P^{opt}_{i,j,k}(x_i\pm\Dx/2)=c_LP_L(x_i\pm\Dx/2)+c_CP_C(x_i\pm\Dx/2)+c_RP_R(x_i\pm\Dx/2).
\eeq
The value of these coefficients, also called ``optimal weights" is $c_L=c_R=\frac{1}{6}$, $c_C=\frac{2}{3}$ for the fourth-order CWENO scheme \cite{LPR99}.

For non-smooth problems, the idea is to replace these weights $c_L, c_C$ and $c_R$ by other positive weights $w_L, w_C$ and $w_R$ (with $w_L+w_C+w_R=1$):
\beq
\label{eq:pr}
R_{i,j,k}(x_i\pm\Dx/2)=w_LP_L(x_i\pm\Dx/2)+w_CP_C(x_i\pm\Dx/2)+w_RP_R(x_i\pm\Dx/2),
\eeq
with the weights such that in smooth regions $(w_{\vfa})=(c_{\vfa}), \vfa \in \{L,C,R\}$ so as to have the reconstruction polynomial $R=P^{opt}$, providing high order of accuracy, but in regions with strong discontinuities, the polynomial(s) in which the discontinuity is present is/are associated with a vanishing weight. For example, in Fig. \ref{fig:whyweno}, only $P_R$ would have a non-negligible weight. This can be done through:
\beq
\label{eq:w}
w_{\vfa}=\frac{\alpha_{\vfa}}{\sum \alpha_j}\:\:\:\mathrm{with}\:\:\:\alpha_{\vfa}=\frac{c_{\vfa}}{(\epsilon+IS_{\vfa})^p},
\eeq
where $\vfa \in\{L,C,R\}$, $p=2$ as in \cite{JIS96,LPR99}, $\epsilon$ is a small positive number taken in order to avoid that the denominator becomes zero, typically $\epsilon=10^{-6}$, and $IS_{\vfa}$ corresponds to a measure of the smoothness of the polynomial $P_{\vfa}, \vfa \in\{L,C,R\}$:
\beq
IS_{\vfa}=\sum_{l=1}^{n} \int_{x_{i-1/2}}^{x_{i+1/2}} \Dx^{2l-1} \left( \frac{\partial^l P_{\vfa}}{\partial x^l} \right)^2 \dx.
\eeq
This smoothness indicator is a measure of the total variation of the polynomial, considering as well the higher order variations, while the $\Dx^{2l-1}$ term is there to remove the $\Dx$ terms coming from the polynomial's derivatives. All algebra done, this gives, with the subscripts $(j,k)$ implicitly assumed:
\beqa
\label{eq:isl}
IS_L&=&\frac{13}{12}(Q_{i-2}-2Q_{i-1}+Q_{i})^2+\frac{1}{4}(Q_{i-2}-4Q_{i-1}+3Q_{i})^2,\\
IS_C&=&\frac{13}{12}(Q_{i-1}-2Q_{i}+Q_{i+1})^2+\frac{1}{4}(Q_{i-1}-Q_{i+1})^2,\\
\label{eq:isr}
IS_R&=&\frac{13}{12}(Q_{i}-2Q_{i+1}+Q_{i+2})^2+\frac{1}{4}(3Q_{i}-4Q_{i+1}+Q_{i+2})^2.
\eeqa
And the reconstructed values:
\beqa
\label{eq:Qw}
Q^{\As,W}_{i+1/2}&=&R_{i}(x_i+\Dx/2),\\
&=&\frac{1}{6}[w_L(2Q_{i-2}-7Q_{i-1}+11Q_{i})+w_C(-Q_{i-1}+5Q_{i}+2Q_{i+1})\nonumber\\
&+&w_R(2Q_{i}+5Q_{i+1}-Q_{i+2})],\nonumber\\
\label{eq:Qe}
Q^{\As,E}_{i-1/2}&=&R_{i}(x_i-\Dx/2),\\
&=&\frac{1}{6}[w_L(-Q_{i-2}+5Q_{i-1}+2Q_{i})+w_C(2Q_{i-1}+5Q_{i}-Q_{i+1})\nonumber\\
&+&w_R(11Q_{i}-7Q_{i+1}+2Q_{i+2})].\nonumber
\eeqa
Here the superscripts $W$ and $E$ stand respectively for \lq\lq West\rq\rq\ and \lq\lq East\rq\rq, which correspond to states obtained for $x \to x^{-}_{i+1/2}$ and $x \to x^{+}_{i-1/2}$ respectively.
Please note that other choices for the smoothness indicators (for example \cite{ZHS07,BCC08}) and/or the weights (see for example \cite{HAP05}) are possible, which can reduce the dissipation, improve the accuracy near smooth extrema and/or provide a better convergence to steady states.

When considering a system of equations, it is furthermore important that all dependent variables \lq\lq feel\rq\rq\ the discontinuities and shocks at the same locations. If each variable is reconstructed independently, oscillations can still occur. In order to avoid them, it is advisable to use so-called \lq\lq Global Smoothness Indicators\rq\rq\ ($GSI$), by taking an appropriate combination of the individual smoothness indicators of each variable, and using this combination to compute the weights $w_{\vfa}$ which are common for all variables. For example, for a pure hydrodynamical problem (with $\vB=\vec{0}$), it has been shown that using the smoothness indicator of the density alone for all variables is a good choice \cite{LPR99}, and for MHD problems using a normalised average of the smoothness indicators from the density and the magnetic field components gives good results \cite{BAT06,VAL18}. An illustration of this fact is also shown in section \ref{sec:BrioWu}.
The CWENO reconstruction procedure can hence be computed in practice through the following steps:
\begin{enumerate}
\item use \eqref{eq:isl}-\eqref{eq:isr} to compute individual smoothness indicators for a well chosen ensemble of variables,
\item combine them appropriately to compute global smoothness indicators $(GSI_m), m \in \{L,C,R\}$,
\item compute the weights using \eqref{eq:w} with the $(GSI_m)$,
\item apply these common weights in \eqref{eq:Qw}-\eqref{eq:Qe} to reconstruct all the variables.
\end{enumerate}

\subsection{Reconstruction of the magnetic field components}\label{sec:Bcell}
As mentioned in section \ref{sec:CT}, in the CT approach, each magnetic field component is known primarily as an area average normal to its respective direction (see Fig. \ref{fig:vardef}). In order to determine the fluxes however, each component has to be estimated on the other faces as well. In order to do this, one can first compute the volume-averaged magnetic field inside the cells and deduce the components on any other face through the reconstruction method described above in section \ref{sec:WENO}. The easiest way to determine this volume average is through a polynomial interpolation. Namely, for a $2n^{th}$ order method, one can find the unique polynomial $P_{B,\As}$ of degree at most $2n-1$ which verifies
\beq
P_{B,{\As}}(x_{i+\vfa-1/2})=(B^{\As}_x)_{i+\vfa-1/2,j,k}, \forall \vfa \in [ -n+1, n ]
\eeq
Then, the integral $\frac{1}{\Dx}\int_{x_i-1/2}^{x_i+1/2} P_{B,{\As}}(x) \dx$ gives a $2n^{th}$ order approximation of $\frac{1}{\Dx\Dy\Dz}\int_{\Omega_{i,j,k}} B_x \dx$ 
All derivations done, for a fourth-order scheme, this gives:
\beqa
\frac{1}{\Dx\Dy\Dz}\int_{\Omega_{i,j,k}} B_x \dx&=&\frac{1}{24}[-(B^{\As}_x)_{i-3/2,j,k}-(B^{\As}_x)_{i+3/2,j,k} \nonumber \\
&+&13\left( (B^{\As}_x)_{i-1/2,j,k}+(B^{\As}_x)_{i+1/2,j,k})\right)].
\eeqa
The same method can of course be applied to obtain the volume averages of the other magnetic field components.
Please note that the use of a simple interpolation is possible here since the area-averaged magnetic field components are continuous along their respective directions \cite{LOZ00}. In order to enhance the stability of the scheme, one could however, similarly to the CWENO volume-to-area average reconstructions, use a non-oscillatory reconstruction, as in reference \cite{BMD13}.

When implementing such a scheme in practice, please note that the $B_{\vfa}$ component of the magnetic field does not have to be reconstructed along the $\vfa$ direction since it is known as a primary data. This also means that $(B^{\As,E}_x)_{i\pm 1/2,j,k}=(B^{\As,W}_x)_{i\pm 1/2,j,k}=(B^{\As}_x)_{i\pm 1/2,j,k}$ and similarly for the other directions, which is consistent with the continuity of $B_m$ along direction $m$.
\subsection{Solving the Riemann problem}
\label{sec:RP}
Through the reconstruction procedure, two states for $(\vU,\vB)$ are determined at each interface between cells: a ``West" and an ``East" one. At each interface, we are thus facing a so-called \lq\lq Riemann problem\rq\rq\footnote{A Riemann problem is an initial value problem with a system in a certain constant state on one side of an impermeable interface and a different constant state on the other side, presenting hence a discontinuity.
The central question is the system's temporal evolution after removal of the interface.} and want to estimate the fluxes through this interface. We do not detail here how such problems are solved in practice, but the interested reader can see for example \cite{TOR09} for details. Numerically, several Riemann solvers with different levels of complexity have been proposed, see for example the references \cite{RUS61,ROE81,KNP01,MIK05}.
We present here the simplest one, the so-called \lq\lq Rusanov flux\rq\rq\ \cite{RUS61} (also known as \lq\lq Local Lax-Friedrichs\rq\rq\ or LLF \cite{LOZ00})
 
\beq
\label{eq:LLFflux}
\vf^x_{LLF}=\frac{1}{2}\left(\vf^x(\vU^E,\vB^E)+\vf^x(\vU^W,\vB^W)\right)-\frac{a^x}{2}(\vU^E-\vU^W),
\eeq
where $\vf^x$ is the physical flux in the $\vx$-direction (see \eqref{eq:fluxfx}) and $a^x$ is the maximum local speed of propagation of information in the system in the $\vx$-direction, which corresponds in the case of the MHD equations to the fast-magnetosonic wave:
\beq
\label{eq:an}
a^x=\max\left( (|v_x|+c^x_f)^W,(|v_x|+c^x_f)^E\right),
\eeq
with $c^x_f$ the magneto-sonic speed:
\beq
c^x_f=\sqrt{\frac{1}{2}\left( (c_s^2+c_A^2)+\sqrt{(c_s^2+c_A^2)^2-4c_s^2\frac{B_x^2}{\rho}} \right)},
\eeq
where $c_s=(\gamma p/\rho)^{1/2}$ and $c_A=(|\vb|^2/\rho)^{1/2}$ are the sound speed and the Alfv\'{e}n speed, respectively. The Rusanov flux has the advantage of being very simple and computationally inexpensive. However, it is way more dissipative than other more accurate Riemann solvers, such as the Roe \cite{ROE81} or the HLLD solver \cite{MIK05}. 
Please note: if one wants to have a scheme of order higher than two, one cannot simply use \eqref{eq:LLFflux} with the computed area averages. The next section addresses this issue.
\subsection{Passage through point values}
\label{sec:pointvalues}
As can be seen from \eqref{eq:fluxfx}, quantities such as $\iint \rho v_x^2$ need to be known in order to compute the fluxes. From the reconstruction step, the area averages $\iint \rho$ and $\iint \rho v_x$ are known. Mathematically however, $\iint \rho v_x^2 \neq \frac{(\iint \rho v_x)^2}{\iint \rho}$. A way to walk-around this difficulty is to compute point values from the area averages, since, for point values, $\rho v_x^2=\frac{(\rho v_x)^2}{\rho}$. The main idea to perform such a transformation comes from a Taylor expansion and can be found in \cite{COC11,BUH14}. For the sake of simplicity, we derive first the formula for a one-dimensional problem:
\beqa
\label{eq:TayFtoP}
Q_{i}&=&\frac{1}{\Dx}\int_{\epsilon=-\Dx/2}^{\Dx/2} q(x_i+\epsilon) \deps,\\
\label{eq:TayFtoP-2}
&=&\frac{1}{\Dx}\int \left(q(x_i) +\epsilon q'(x_i)+\frac{\epsilon^2}{2!} q''(x_i)+\frac{\epsilon^3}{3!}q'''(x_i)+O(\epsilon^4) \right) \deps,\\
\label{eq:TayFtoP-3}
&=&q(x_i) +\frac{\Dx^2}{24} q''(x_i)+O(\Dx^4).
\eeqa
Through this formula, one can also see that $Q_i=q_i+O(\Dx^2)$. Thus, identifying the spatial average with a point value in the middle of the considered domain results in an error of second order. This is why we need to consider the next term, $\frac{\Dx^2}{24} q''(x_i)$ if we want to have a fourth-order scheme. It can be shown 
that a second-order approximation of $q''(x_i)$ is \cite{BUH14}:
\beq
\label{eq:2ndderiv}
q''(x_i)=\frac{q_{i-1}-2q_i+q_{i+1}}{\Dx^2}+O(\Dx^2)=\frac{Q_{i-1}-2Q_i+Q_{i+1}}{\Dx^2}+O(\Dx^2),
\eeq
where the first equality contains point values and the second cell averages. Please note that, except in this
special case, one cannot in general simply replace point values by cell averages in the formulas.

The Taylor expansion of \eqref{eq:TayFtoP-2} can be generalized for more than one dimension. For a 2D fourth-order area-average-to-point transformation, this leads to the following expression:
\beqa
\label{eq:FtoP}
q_{i \pm 1/2,j,k}&=&Q_{i \pm 1/2,j,k} - \frac{Q_{i \pm 1/2,j+1,k}-2Q_{i \pm 1/2,j,k}+Q_{i \pm 1/2,j-1,k}}{24}\\
\nonumber	&-& \frac{Q_{i \pm 1/2,j,k+1}-2Q_{i \pm 1/2,j,k}+Q_{i \pm 1/2,j,k-1}}{24} +O(\Dx^4\!+\!\Dy^4\!+\!\Dz^4).
\eeqa
Please note that for schemes of order higher than 4, cross derivatives (for example $\partial^2 x \partial^2 y$) need to be considered. 

These computed point values can be used in order to determine point-valued fluxes through \eqref{eq:LLFflux}. These point-valued fluxes can then be transformed using a point-to-area-average transformation, obtained by using \eqref{eq:2ndderiv} and \eqref{eq:TayFtoP-3}:
\beqa
\label{eq:PtoF}
Q_{i \pm 1/2,j,k}&=&q_{i \pm 1/2,j,k} + \frac{q_{i \pm 1/2,j+1,k}-2q_{i \pm 1/2,j,k}+q_{i \pm 1/2,j-1,k}}{24}\\
\nonumber	&+& \frac{q_{i \pm 1/2,j,k+1}-2q_{i \pm 1/2,j,k}+q_{i \pm 1/2,j,k-1}}{24} +O(\Dx^4\!+\!\Dy^4\!+\!\Dz^4).
\eeqa
Plugging $\vf^x_{LLF}$ obtained through \eqref{eq:LLFflux} in this relation gives a higher-order approximation to the area-averaged flux $\vF^{\As,x}$, which is used in \eqref{eq:patvU} in order to evolve the cell averages of the hydrodynamical quantities in time (see section \ref{sec:timeINT} for the time integration).

Once again, please note that the fact that \eqref{eq:PtoF} and \eqref{eq:FtoP} look very similar to each other is a special case, valid only because of the similarity in relation \eqref{eq:2ndderiv}. This similarity is lost for schemes of order higher than 4.
\subsection{Electric fluxes on the edges}
\label{sec:Ee}
\begin{figure}
\begin{center}
\includegraphics[width=1.\textwidth]{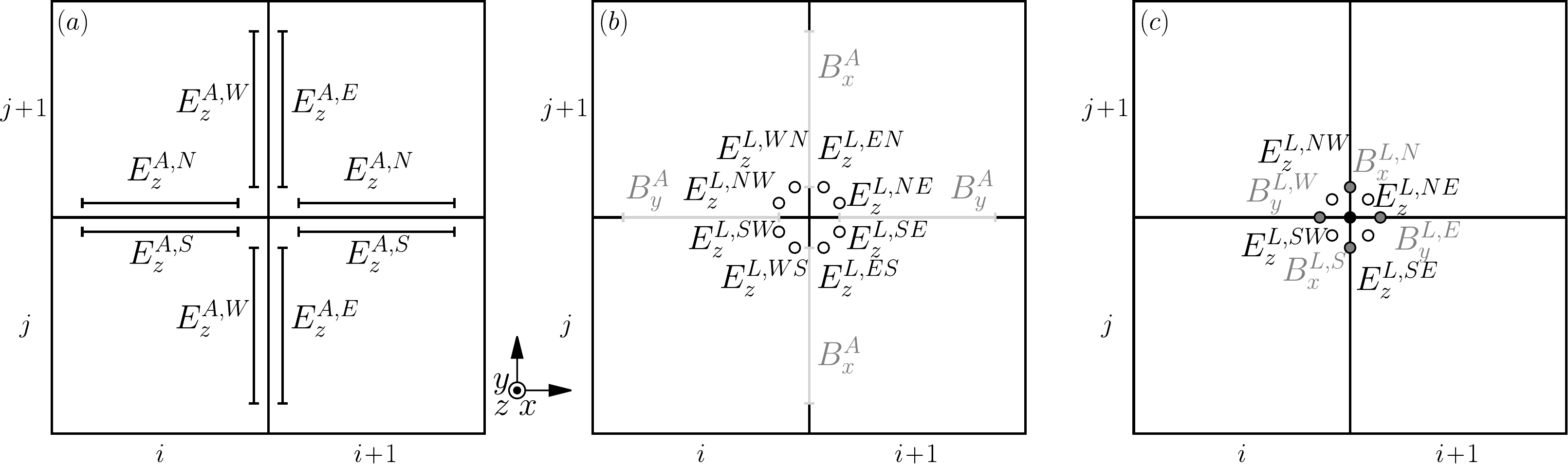}
\caption{ \label{fig:eecomp} Illustration of the computation of $(E^{\Ls}_z)_{i+1/2,j+1/2,k}$ (black circle in subgraph $(c)$): $(a)$ CWENO line-average reconstruction of the area-averaged electric field, $(b)$ CWENO line-average reconstruction of the area-averaged magnetic field and taking the mean of the two possibilities for the electric field in each quadrant, $(c)$ Deducing the electric flux thanks to the multidimensional Riemann solver.}
\end{center}
\end{figure}
In order to evaluate the right hand side for the magnetic field evolution, we have to take care of a two-dimensional Riemann problem. Indeed, formally, \eqref{eq:mag} can be reformulated in conservative form:
\begin{equation}
\pat \vb= - \vnabla \cdot \begin{pmatrix} 0 & -E_z & E_y \\ E_z & 0 & -E_x \\ -E_y & E_x & 0 \end{pmatrix}.
\end{equation}
This means that the electric field plays the role of a flux. Thanks to section \ref{sec:pointvalues}, electric field point values can be computed from the velocity and magnetic field point values and then transformed into area averages, two per interface between two cells, for example $\vE^{\As,W}_{i+1/2,j,k}$ and $\vE^{\As,E}_{i+1/2,j,k}$ at the interface between cells $(i,j,k)$ and $(i+1,j,k)$. Then, the same CWENO procedure as the one described in section \ref{sec:WENO} can be applied to these area averages in the two remaining directions in order to deduce line averages: in our example we would obtain for a reconstruction in the $\vy$-direction the four states $\vE^{\Ls,SW}_{i+1/2,j+1/2,k}, \vE^{\Ls,SE}_{i+1/2,j+1/2,k}, \vE^{\Ls,NW}_{i+1/2,j-1/2,k}$ and $\vE^{\Ls,NE}_{i+1/2,j-1/2,k}$ where the superscripts $S$ and $N$ stand respectively for ``South" ($y\to y_{j+1/2}^-)$ and ``North" ($y\to y_{j-1/2}^+$). For these area-to-line CWENO reconstructions, one can take for \eqref{eq:Qw}-\eqref{eq:Qe} the weights computed by taking the mean of the $GSI$ in the reconstruction direction of the two cells surrounding the area (see section \ref{sec:WENO}). In the end, at each edge surrounded by four faces, we have four different states for the electric field (see Fig. \ref{fig:eecomp}), which are computed by taking the mean of the two possible reconstructed values (for example, one can obtain $\vE^{\Ls,SW}_{i+1/2,j+1/2,k}$ by $(a)$ reconstructing $\vE^{\As,W}_{i+1/2,j,k}$ along $\vy$ or $(b)$ reconstructing $\vE^{\As,S}_{i,j+1/2,k}$ along $\vx$). These four states are then plugged in a 2D Riemann solver. We present here the LLF variant of the multidimensional Riemann solver of reference \cite{BAL10}:
\beqa
(E^{\Ls}_z)&=& \frac{1}{4}\Big[(E^{\Ls}_z)^{SE}+(E^{\Ls}_z)^{SW}+(E^{\Ls}_z)^{NE}+(E^{\Ls}_z)^{NW}\Big]\\
\nonumber
&&+\frac{S^a}{2}\left({(B_y^{\Ls})^E-(B_y^{\Ls})^W}\right)-\frac{S^a}{2}\left({(B_x^{\Ls})^N-(B_x^{\Ls})^S}\right),
\label{eq:eez}
\eeqa
where the quantities with a superscript $\Ls$ are computed through area-to-line average reconstructions and are indicated in Fig. \ref{fig:eecomp} and $S^a$ is the maximum speed of propagation of information, estimated as:
\beq
S^a_{i+1/2,j+1/2,k}=\max(a^x_{i+1/2,j,k},a^x_{i+1/2,j+1,k},a^y_{i,j+1/2,k},a^y_{i+1,j+1/2,k}).
\eeq
This choice for evaluating $S^a$ is done for reasons of efficiency, since the local maximum speeds $a^{x,y}$ have already been computed in order to solve the 1D Riemann problems (see section \ref{sec:RP}), as mentioned in \cite{LOZ04}.
Please note that refinements of this solver are available in the literature: in reference \cite{BAL10}, a solver with four different speeds can be found and further improvements have been made in reference \cite{BNK17}.
\subsection{Summary: the complete procedure to determine the R.H.S}
The numerical method to compute the right hand sides of \eqref{eq:patvU} and \eqref{eq:dbSA}-\eqref{eq:dbSB} can be summarized as follows:
\begin{enumerate}
\item Compute cell averages for the magnetic field components (section \ref{sec:Bcell}).
\item For each dimension $n \in \{x,y,z\}$:
	\begin{enumerate}
		\item reconstruct area averages of all quantities but $B_n$ along $n$ (sections \ref{sec:WENO} and \ref{sec:Bcell}),
		\item transform the area averages to point values (section \ref{sec:pointvalues}),
		\item deduce point-valued fluxes (section \ref{sec:RP}) and electric field (section \ref{sec:Ee}),
		\item transform the point-valued fluxes and electric field to area averages (section \ref{sec:pointvalues}),
		\item deduce one term of \eqref{eq:patvU}.
	\end{enumerate}
\item For each dimension $n \in \{x,y,z\}$:
	\begin{enumerate}
		\item reconstruct line averages of the $\vB$ and $\vE$ components normal to $n$ (section \ref{sec:Ee}),
		\item deduce through \eqref{eq:eez} two terms in \eqref{eq:dbSA}-\eqref{eq:dbSB}.
	\end{enumerate}
\end{enumerate}
Please note that the dimension-by-dimension approach presented here is computationally more efficient than an actual multi-dimensional CWENO reconstruction which requires bi- (resp. tri)-quadratic polynomials in 2D (resp. 3D)\cite{VEM18}. 

\section{Time integration}
\label{sec:timeINT}
Numerous explicit time integration methods have been developed with different mathematical properties making them more appropriate in different contexts. When dealing with problems with strong discontinuities and shocks, the so-called Strong Stability Preserving Runge-Kutta (SSPRK) methods are advantageous since they avoid
additional oscillations generated by the time integration step \cite{GOS98,GST01}. Furthermore, their linear stability region is typically greater than standard Runge-Kutta schemes, allowing larger timesteps and hence a large gain in efficiency.

The numerical tests presented in this chapter are done using the ten-stage fourth-order SSPRK method presented in pseudocode 3 of reference \cite{KET08}, which is repeated for the sake of completeness in algorithm 1.
\begin{algorithm}
\caption{10 stages SSPRK4 method from \protect\cite{KET08}.}
\begin{algorithmic}
\STATE $k_1 \leftarrow w$
\FOR{$s=1:5$}
\STATE	$k_1 \leftarrow k_1+\frac{\Dt}{6} F(k_1)$
\ENDFOR
\STATE $k_2 \leftarrow \frac{1}{25}w+\frac{9}{25}k_1$
\STATE $k_1 \leftarrow 15k_2-5k_1$
\FOR{$s=6:9$}
\STATE	$k_1 \leftarrow k_1+\frac{\Dt}{6} F(k_1)$
\ENDFOR
\STATE $w \leftarrow k_2+\frac{3}{5}k_1+\frac{\Dt}{10}F(k_1)$
\end{algorithmic}
\end{algorithm}
In this pseudocode, $\Dt$ is the timestep, $w$ represents the whole state of the system (the hydrodynamical variables as well as the magnetic field) initially at instant $t$ but at the last line at instant $t+\Dt$, $F(k_1)$ is the right-hand side computed from the state $k_1$ and $s$ the stage number. Three registers are needed: $k_1, k_2$ and one to store $F(k_1)$. In order to obtain a stable solution with correct results, the timestep $\Dt$ cannot be arbitrarily large but is restricted by the so-called Courant-Friedrichs-Lewy (CFL) stability criterion \cite{CFL28}. We present here a formulation of this criterion that can be found in \cite{TIT05} for example:
\beq
\label{eq:cfl}
\Dt \leq C_{CFL} \min_{i,j,k}(\frac{\Dx}{a^x_{i,j,k}},\frac{\Dy}{a^y_{i,j,k}},\frac{\Dz}{a^z_{i,j,k}}),
\eeq
where $C_{CFL}$ is a constant called the Courant number, $a^n (n\in\{x,y,z\})$ is the maximum speed of propagation of information in the $\vn$-direction, which in our case corresponds to fast magnetosonic waves and is expressed in \eqref{eq:an}. In the case of the ten-stage method used here, even though a high number of stages are performed, its stability region is very large, allowing a high Courant constant $C_{CFL}=1.95$ for the 2D and $C_{CFL}=1.55$ for the 3D systems we consider in this chapter.

In order to spare computer time, it is advisable to allow the timestep to vary. At each iteration, the maximum allowed timestep verifying \eqref{eq:cfl} is computed during the first stage of the SSPRK method and kept the same for the remaining stages, until the next iteration.
\section{Strong shocks and negative pressure/density}
\label{sec:fallback}
Even when using a non-oscillatory reconstruction and a SSPRK time-integration method, negative densities or pressure may still appear when dealing with strong discontinuities or shocks. The increase of the stiffness of the reconstruction polynomials together with the reduction of numerical dissipation which are both consequences of the increase of the algorithm's precision make high order schemes particularly vulnerable to this issue. Indeed, even though numerical dissipation causes amplitude errors in the solution, which tend to smooth out small-scale structures, it has the positive side effect of smoothening out unphysical reconstruction variations as well. Because of this, lower order schemes deliver less precise solutions but are more robust when compared to higher order schemes.
A negative density or pressure can arise after the reconstruction step (sec. \ref{sec:WENO}), the area-to-point transformation (sec. \ref{sec:pointvalues}) or the time integration (sec. \ref{sec:timeINT}).

A possibility to deal with this issue is to use a fallback approach, which consists of using lower-order reconstruction procedures in the vicinity of strong shocks, so as to smooth out the discontinuity and avoid additional oscillations.

First of all, one can use an \textit{a-posteriori} fallback approach: if after the reconstruction procedure (section \ref{sec:WENO}) a negative density or pressure is obtained, then a lower-order reconstruction is performed locally, such as a second-order TVD limiter \cite{VLE77}, or in the extreme case Godunov's scheme $\vU^{\As}_{i\pm 1/2,j,k}=\vU_{i,j,k}$. Similarly, if the area-to-point transformation would give an unphysical result, it can be switched off locally: the area averages are identified with the point values.

However, this may not be enough for a strongly shocked problem. An idea then is to use an \textit{a-priori}
fallback approach, also called ``flattening" \cite{COW84,BAL12}, which would reconstruct at lower order even though the higher order methods still give a physical state such as to avoid the appearance of unphysical states in the next timesteps. In practice, assume a high-order reconstruction method $HO$ and a low-order reconstruction method $LO$, which give the reconstructed states $(\vU^{\As},\vB^{\As})^{HO}_{i\pm 1/2,j,k}$ and $(\vU^{\As},\vB^{\As})^{LO}_{i\pm 1/2,j,k}$ respectively. Then, the reconstructed state at the interfaces $(i \pm 1/2,j,k)$ is given by:
\begin{equation}
(\vU^{\As},\vB^{\As})_{i\pm 1/2,j,k}=w^f_{i,j,k}(\vU^{\As},\vB^{\As})^{HO}_{i\pm 1/2,j,k}+(1-w^f_{i,j,k})(\vU^{\As},\vB^{\As})^{LO}_{i\pm 1/2,j,k},
\end{equation}
with a weight $w^f_{i,j,k} \in [0,1]$ called the \lq\lq flattener\rq\rq, equal to 1 in smooth regions and going to 0 as the considered region is more and more shocked. In a similar way, this flattening also occurs for the area-to-point transformation in case the lower-order method is first- or second-order accurate:
\beqa
\label{eq:FTOPflat}
\nonumber
q_{i \pm 1/2,j,k}& \approx &Q_{i \pm 1/2,j,k} - w^f_{i\pm 1/2,j,k}\frac{Q_{i \pm 1/2,j+1,k}-2Q_{i \pm 1/2,j,k}+Q_{i \pm 1/2,j-1,k}}{24}\\
	&-& w^f_{i\pm 1/2,j,k}\frac{Q_{i \pm 1/2,j,k+1}-2Q_{i \pm 1/2,j,k}+Q_{i \pm 1/2,j,k-1}}{24},
\eeqa
and similarly for the point-to-area transformation. There are several possible choices for the flattener: in the XTROEM-FV code \cite{ROM16}, the Jameson indicator \cite{JST81} is used. Another flattener, well suited for the structure of the MHD equations, can be found in \cite{BAL12}, which also presents refinements in order to include cells that are not affected by the shock yet but will shortly be. We present however here one inspired by reference \cite{COW84} which uses the pressure gradient. Namely, we define first, for a reconstruction in the $\vx$-direction for example, the shock indicator:
\beq
	s^x_{i,j,k}=\frac{|\tilde{p}_{i+1,j,k}-\tilde{p}_{i-1,j,k}|}{\tilde{p}_{i,j,k}},
\eeq
with $\tilde{p}_{i,j,k}$ an estimate of the pressure in cell $(i,j,k)$ computed by considering the cell averages as if they were point values. Then, we choose two thresholds $\tau^{HO} < \tau^{LO}$ and define the flattener $w_f$ as:
\beq
w^f_{i,j,k} = \begin{cases} 1, & \mbox{if } s^x_{i,j,k}<\tau^{HO}, \\ 1-\frac{s^x_{i,j,k}-\tau^{HO}}{\tau^{LO}-\tau^{HO}}, & \mbox{if } \tau^{HO} \leq s^x_{i,j,k} \leq \tau^{LO}, \\ 0, & \mbox{if } \tau^{LO}<s^x_{i,j,k}. \end{cases}
\eeq
The choice of the thresholds $\tau^{HO}$ and $\tau^{LO}$ is dependent both on the problem and the reconstruction method. If the lower-order method is second-order, a reasonable first try is $\tau^{HO},\tau^{LO} \lesssim 4$. Indeed, if one would use naively a linear approximation, without limiters, to reconstruct the pressure at the cell interfaces: $p_{i \pm 1/2,j,k}=\tilde{p}_{i,j,k} \pm \frac{\tilde{p}_{i+1,j,k}-\tilde{p}_{i-1,j,k}}{2 \Dx} \cdot \frac{\Dx}{2}$, then a value of $s^x_{i,j,k}=4$ would already give a zero-pressure on one of the interfaces. The distance between $\tau^{LO}$ and $\tau^{HO}$ is there to provide a smooth transition from high to low order. In practice, this rough upper boundary may be too high for highly shocked problem and may need to be decreased severely depending on the dynamical state of the system.

For the flattener at the area-to-point-values transformations in eq. \ref{eq:FTOPflat}, one can take a function of the $w^f$ used for reconstruction of the neighbouring cells. Adding a superscript $n \in \{x,y,z\}$ corresponding to the reconstruction direction, we decided to use:
\beq
	w^f_{i+1/2,j,k}=\min(w^{f,y}_{i,j,k},w^{f,y}_{i+1,j,k},w^{f,z}_{i,j,k},w^{f,z}_{i+1,j,k})
\eeq
This method is of course generalisable to more than two reconstruction methods and area$\leftrightarrow$point values transformations. Another possibility to deal with strong discontinuities may be WENO schemes which automatically adapt their order in their vicinity \cite{BGS16}.

\section{Numerical tests}
\label{sec:numtests}
In the case of the numerical scheme presented here, the following aspects have to be verified:
\begin{enumerate}
\item the scheme should converge towards an exact solution in continuous space and time at the proper rate
 in the asymptotic limit of  $\Delta x$, $\Delta t\rightarrow 0$, $\Delta x/\Delta t \neq 0$. This is the object of section \ref{sec:ctest}.
\item The solenoidality of the magnetic field should be preserved up to machine precision: in all the tests performed here, $\nabla \cdot \vB$ is close to machine precision and does not grow over time. Since the $\nabla \cdot \vB$ is \textit{preserved} in the CT approach, we would like to point out that the initialization is very important in this respect. In practice, it is a good idea to use analytical expressions for the area averages of the magnetic field components, and not point values in the middle of the faces for example, even when using a second-order scheme.
\item The scheme should be robust, that is, be able to handle correctly problems with shocks and discontinuities. In particular, it should not give unphysical values such as a negative pressure or density. This is the object of section \ref{sec:numtests_shocks}.
\end{enumerate}

In all the tests, we take $\gamma=\frac{5}{3}$, which corresponds physically to a mono-atomic gas with three degrees of freedom,
and $C_{CFL}=1.95$ for 1D and 2D problems and $C_{CFL}=1.55$ for 3D problems, unless mentioned otherwise. The boundary conditions are periodic in all tests, apart for the Brio-Wu Riemann problem.


\subsection{Verification of the scheme's order}
\label{sec:ctest}
For a smooth problem, the error between the solution given by the numerical scheme and the exact solution is expected to diminish asymptotically at a well-defined rate corresponding to the scheme's order. A possibility to perform a convergence test is to take a smooth problem with a periodic motion and compare the solution after a period $t=T$ to the initial conditions at $t=0$. When the exact solution is not known, the convergence order can be computed by comparing to a very high resolution run which is then taken as the reference, see the convergence test presented in section \ref{sec:OT}. There are several possibilities to compute the discretization error. We present here the same one as in \cite{MTB10} by using the $L_1$ norm:
\beq
\label{eq:ctest}
\delta \vW=\frac{1}{N_xN_yN_z}\sum_{i,j,k} |\vW_{i,j,k}(t=T)-\vW_{i,j,k}(t=0)|
\eeq
with $\vW=(\vU,\vB)$ so that $\delta \vW$ is the error vector for all the eight MHD variables. $N_x, N_y$ and $N_z$ are the number of gridpoints in each direction. We define furthermore $\delta W_{mean}=\frac{1}{8}\sum_i \delta \vW_i$ the mean of these errors. In order to measure the convergence, the same numerical test is repeated for several resolutions $r_1, r_2, ..., r_N$ and the experimental order of convergence $EOC$ is given by:
\beq
EOC_{\vfa}=\frac{|\log(\delta W_{mean}(r_{\vfa}))-\log(\delta W_{mean}(r_{\vfa-1}))|}{|\log(r_{\vfa})-\log(r_{\vfa-1})|}
\eeq
and should be asymptotically equal to the scheme's order $n$.

\subsection{Smooth problems}
Smooth problems are useful to determine the convergence order of a numerical scheme. We present here a circularly polarized Alfv\'{e}n wave and a 3D MHD vortex. Furthermore, the fourth-order scheme, denoted in the following by CWENO4, is compared with a second-order scheme denoted TVD2 in order to show some advantages of using higher-order numerics. For the second-order scheme, we make use of the TVD limiter of van Leer \cite{VLE77} instead of the fourth-order CWENO reconstruction:

\beqa
Q^{\As,W}_{i+1/2}&=&Q_i+\frac{\max [(Q_{i+1}-Q_i)(Q_i-Q_{i-1}),0]}{Q_{i+1}-Q_{i-1}},\\
Q^{\As,E}_{i-1/2}&=&Q_i-\frac{\max [(Q_{i+1}-Q_i)(Q_i-Q_{i-1}),0]}{Q_{i+1}-Q_{i-1}}.
\eeqa

In addition, the area averages$\leftrightarrow$point values transformations are switched off (that is, the transformations of section \ref{sec:pointvalues} are identity functions). We maintain however the fourth-order magnetic field area-to-volume interpolation as well as the fourth-order time integration so as to be consistent with the a-priori fallback approach for the shocked problems, see section \ref{sec:numtests_shocks}. In the tests performed here, 
the present CWENO4 scheme is about $75\%$ more expensive than
the TVD2 one. This estimate is of course system, optimization and implementation dependent.
The moderate increase of computational effort comes together with a rise in precision
that permits to lower numerical resolution as compared to lower-order techniques.

For these smooth problems, the fallback approach is not needed at all and has been deactivated. 
To show the need of the passage through point values, we consider a CWENO4A scheme, which is the CWENO4 one when ignoring the passage through point values, plugging directly (A)rea averages in the Riemann solver.

\subsubsection{Circularly polarized Alfv\'{e}n wave}
\label{sec:Alfvenwave}

\begin{figure}
\centering
\includegraphics[width=0.9\textwidth]{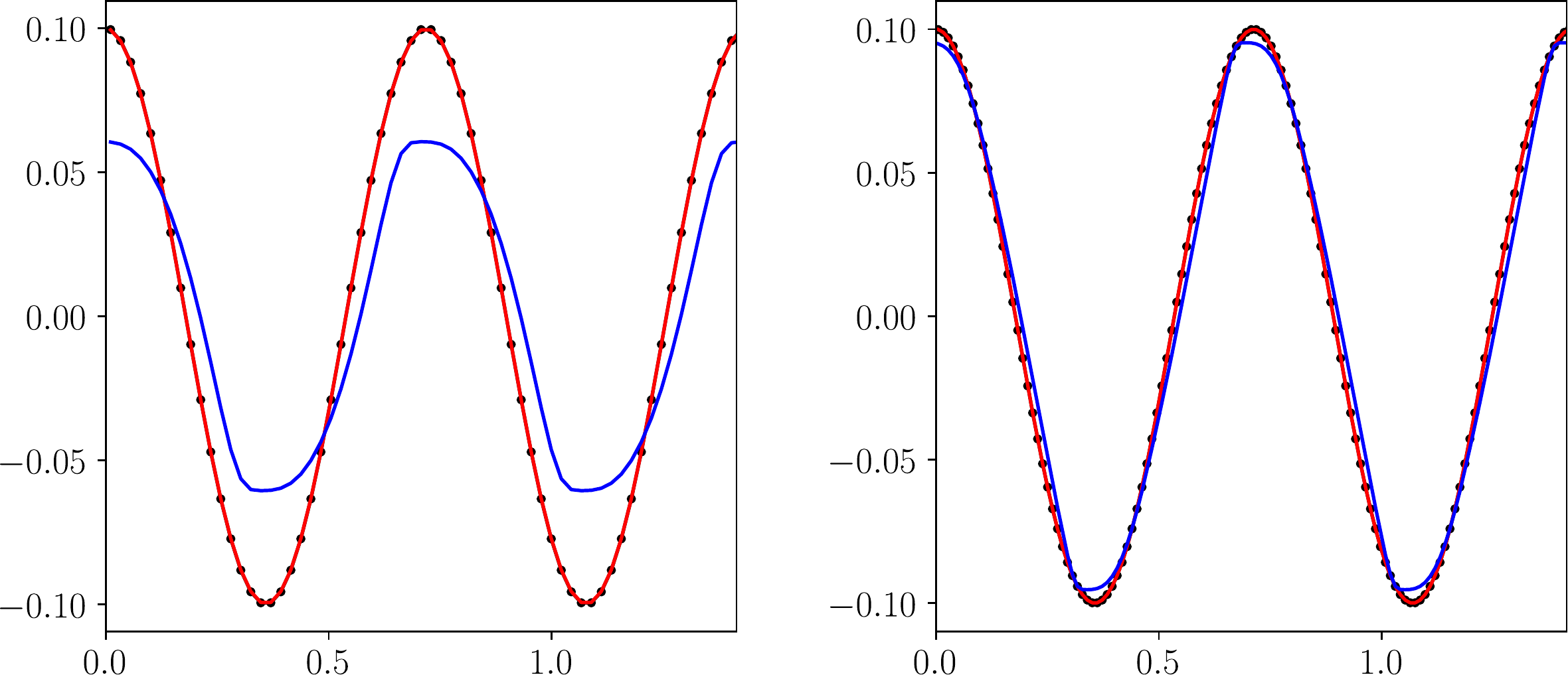}

\begin{tabular}{ cccccccc }
 \hline
\multicolumn{2}{c}{resolution} & $32^2$ & $64^2$ & $128^2$ & $256^2$ & $512^2$ & $1024^2$ \\
\hline
\multirow{3}{*}{TVD2} & $\delta U_{mean}$ & 2.46 $10^{-3}$ & 7.03 $10^{-4}$ & 1.76 $10^{-4}$ & 4.18 $10^{-5}$ & 9.88 $10^{-6}$& 2.31 $10^{-6}$ \\
                        & EOC      & -               & 1.81             & 2.00 & 2.07 & 2.08 & 2.09 \\
                        & $E_{loss}$      & 4.20 $10^{-2}$               & 5.14 $10^{-3}$             & 6.20 $10^{-4}$ & 7.53 $10^{-5}$ & 9.21 $10^{-6}$ & 1.13 $10^{-6}$ \\
\hline
\multirow{3}{*}{CWENO4} & $\delta U_{mean}$ & 1.50 $10^{-5}$ & 4.86 $10^{-7}$ & 1.53 $10^{-8}$ & 5.08 $10^{-10}$ & 1.87 $10^{-11}$& 9.95 $10^{-13}$ \\
                        & EOC      & -               & 4.95             & 4.99 & 4.91 & 4.76 & 4.23 \\
                        & $E_{loss}$      & 7.01 $10^{-4}$               & 2.46 $10^{-5}$             & 7.75 $10^{-7}$ & 2.27 $10^{-8}$ & 5.91 $10^{-10}$ & 1.37 $10^{-11}$ \\
\hline
\multirow{3}{*}{CWENO4A} & $\delta U_{mean}$ & 1.51 $10^{-5}$ & 4.87 $10^{-7}$ & 1.59 $10^{-8}$ & 5.79 $10^{-10}$ & 2.49 $10^{-11}$& 1.42 $10^{-12}$ \\
			 & EOC      & -               & 4.95             & 4.94 & 4.78 & 4.54 & 4.14 \\
                        & $E_{loss}$      & 7.01 $10^{-4}$               & 2.45 $10^{-5}$             & 7.75 $10^{-7}$ & 2.27 $10^{-8}$ & 5.91 $10^{-10}$ & 1.41 $10^{-11}$ \\
\hline
\end{tabular}

\caption{ \label{fig:CPAW} Circularly polarized Alfv\'{e}n wave: cuts of $v_z$ along the main diagonal after 100 periods, (top), left: at resolution $64^2$, right: at resolution $128^2$. The CWENO4 solution (in red) is extremely close to the reference solution at $t=0$ (in black), whereas the TVD2 solution (in blue) presents some significant amplitude and shape error. Bottom: convergence of errors, EOC and numerical dissipation for different schemes after one period.}
\end{figure}

Alfv\'{e}n waves are exact and smooth, inherently linear solutions of the nonlinear MHD equations which are widely used to check the accuracy of a scheme
\cite{TOT00,SGA08,MTB10,ROM16}. We use here the same initial conditions as the circularly polarized Alfv\'{e}n wave in \cite{ZIE04}: $\rho=1$, $\vv=\Big(-\frac{A}{\sqrt{2}}\sin(2\pi(x+y)), \frac{A}{\sqrt{2}}\sin(2\pi(x+y)),A\cos(2\pi(x+y))\Big)$, $\vB=\Big(\frac{B_0}{\sqrt{2}}+\frac{A}{\sqrt{2}}\sin(2\pi(x+y)),\frac{B_0}{\sqrt{2}}-\frac{A}{\sqrt{2}}\sin(2\pi(x+y)),$ $-A\cos(2\pi(x+y))\Big)$ and $p=0.1$, with the wave amplitude $A=0.1$ and the mean-magnetic field $B_0=\sqrt{2}$ in the computational domain $(x,y)\in [0,1]^2$.

The convergence of errors is measured by comparing the solution after one period (at $t=0.5$, since the wave is propagating at the Alfv\'{e}n speed $\sqrt{\frac{B_0^2}{\rho}}=\sqrt{2}$ along the diagonal) to the initial conditions. Figure \ref{tab:CPAW} shows the results for the TVD2, CWENO4 and CWENO4A schemes. The relative energy dissipation $E_{loss}$ is measured by $E_{loss}=\frac{E_{T}(t=0)-E_{T}(t=T)}{E_{T}(t=0)}$ with $E_T$ the total energy (kinetic+magnetic) but removing the mean-fields. We can see that for such a linear solution, even when omitting the passage through point values, the fourth-order scheme gives fourth-order convergence. This is the reason why nonlinear tests, like the one in the next section, have to be performed as well in order to check the convergence order properly. Fig. \ref{fig:CPAW} shows as well some effects of the higher numerical dissipation of the lower order schemes on the wave amplitude after 100 periods, at $t=50$. In order to achieve an as low dissipation as the CWENO4 scheme at resolution $128^2$, the TVD2 scheme needs to go to resolutions higher than $1024^2$, leading to a way higher computational cost, about a factor 300.

\subsubsection{3D MHD vortex}
\label{sec:3DMHV}

\begin{figure}
\centering
\includegraphics[width=0.9\linewidth]{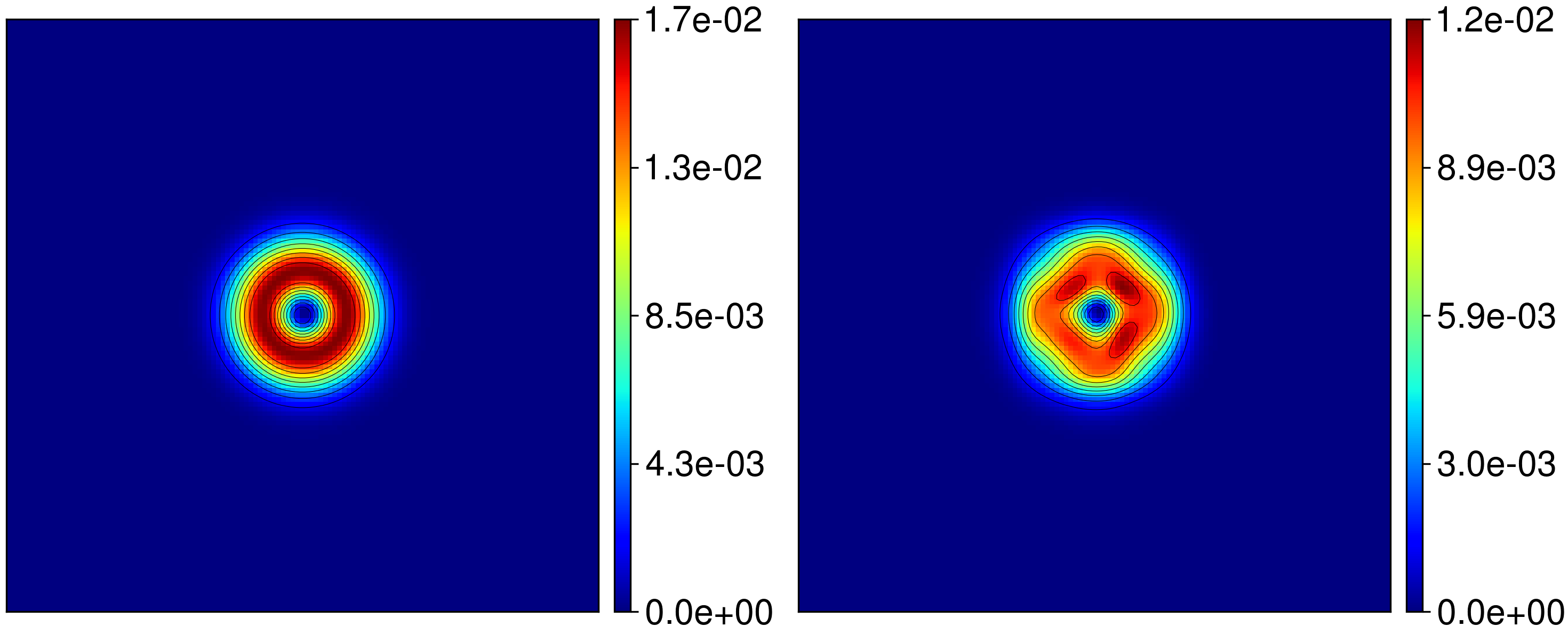}

\begin{tabular}{ cccccccc }
 \hline
\multicolumn{2}{c}{resolution} & $32^3$ & $64^3$ & $128^3$ & $256^3$ & $512^3$ & $1024^3$ \\
\hline
\multirow{3}{*}{TVD} & $\delta U_{mean}$ & 8.85 $10^{-4}$ & 4.33 $10^{-4}$ & 1.35 $10^{-4}$ & 3.92 $10^{-5}$ & 1.02 $10^{-5}$& *** \\
                        & EOC      & -               & 1.03             & 1.68 & 1.79 & 1.94 & *** \\
                        & $E_{loss}$      & 0.92               & 0.58  & 0.15 & 2.02 $10^{-2}$ & 2.39 $10^{-3}$ & *** \\
\hline
\multirow{3}{*}{CWENO4} & $\delta U_{mean}$ & 4.22 $10^{-4}$ & 4.35 $10^{-5}$ & 2.12 $10^{-6}$ & 8.86 $10^{-8}$ & 4.40 $10^{-9}$& *** \\
                        & EOC      & -               & 3.28             & 4.36 & 4.58 & 4.33 & *** \\
                        & $E_{loss}$      & 0.58               & 6.36 $10^{-2}$             & 2.23 $10^{-3}$ & 7.00 $10^{-5}$ & 2.04 $10^{-6}$ & *** \\
\hline
\multirow{3}{*}{CWENO4A}& $\delta U_{mean}$ & 4.21 $10^{-4}$ & 4.37 $10^{-5}$ & 2.30 $10^{-6}$ & 1.64 $10^{-7}$ & 2.98 $10^{-8}$& *** \\
   & EOC      & -               & 3.27             & 4.25 & 3.81 & 2.46 & *** \\
                        & $E_{loss}$      & 0.58                & 6.36 $10^{-2}$             & 2.29 $10^{-3}$ & 7.02 $10^{-5}$ & 2.09 $10^{-6}$ & ***\\
\hline
\end{tabular}
\caption{\label{fig:3DMHVORT} 3D MHD vortex problem: slices of the magnetic pressure at $z\approx 0.04$, resolution $128^3$ after one period (top), left: CWENO4, right: TVD2. Bottom: convergence of errors, EOC and numerical dissipation for different schemes after one period.}
\end{figure}

The MHD vortex problem was first introduced in 2D in \cite{BAL04}. The initial conditions consist of a magnetized vortex structure in force equilibrium that is advected by a velocity field. We present here a 3D extension of this test proposed in \cite{MTB10}. The initial conditions are given by $\rho=1$, $\vv=\Big(1-{y\kappa}\exp\Big[q(1-r^2)\Big],1+{x\kappa}\exp\Big[q(1-r^2)\Big],2 \Big)$, $\vB=\Big({-y\mu}\exp\Big[q(1-r^2)\Big],{x\mu\exp\Big[q(1-r^2)\Big]},0\Big)$ and $p=1+ \frac{1}{4q}\big[\mu^2\big(1-2q(r^2-z^2)\big)-\kappa^2 \rho\big]\exp\Big[2q(1-r^2)\Big]$ with $r^2=x^2+y^2+z^2$, $\kappa=\mu=1/(2\pi)$ and $q=1$ on the computational domain $(x,y,z)\in[-5,5]^3$. Here as well, the solution after one period ($t=10$ since the advection speed is 1 along $\vx$ and $\vy$) is compared to the initial fields. Figure \ref{fig:3DMHVORT} shows that for this nonlinear test problem, omitting the passage through point values leads to a convergence order of 2, which is expected (see eq. \eqref{eq:TayFtoP}). Even though the effect on numerical dissipation is not visible at lower resolutions for CWENO4A, it should become more and more important as the resolution is increased. This figure also shows that in order to obtain a result's quality as the one at resolution $128^3$ for the CWENO4 scheme, the TVD2 scheme needs to be performed at resolution $512^3$, leading to a computational cost increase of about a factor 150.

\subsection{Shocked problems}
\label{sec:numtests_shocks} 
In order to test how well numerical schemes handle shocks, a lot of classical tests have been designed. We present in this chapter the 1D Brio-Wu Riemann problem and the 2D Orszag-Tang vortex. Other classical problems to test the robustness of a numerical scheme include the MHD Rotor Problem (e.g. \cite{BAS99,TOT00,ZIE04,SGA08}),
MHD Blast Waves (e.g. \cite{LOZ00,ZIE04,SGA08,ROM16}), the Current Sheet problem (e.g. \cite{HAS95,ROM16}), the Cloud-Shock Interaction (e.g. \cite{DAW98,TOT00,ZIE04,ROM16}) etc. How an almost identical scheme performs on some of these tests is shown in reference \cite{VAL18}.

Finally, we show an example where the use of an a-priori fallback approach is necessary, namely a supersonic decaying turbulence experiment.

\subsubsection{1D Brio-Wu Riemann problem}
\label{sec:BrioWu}

\begin{figure}
\begin{center}
\includegraphics[width=0.95\textwidth]{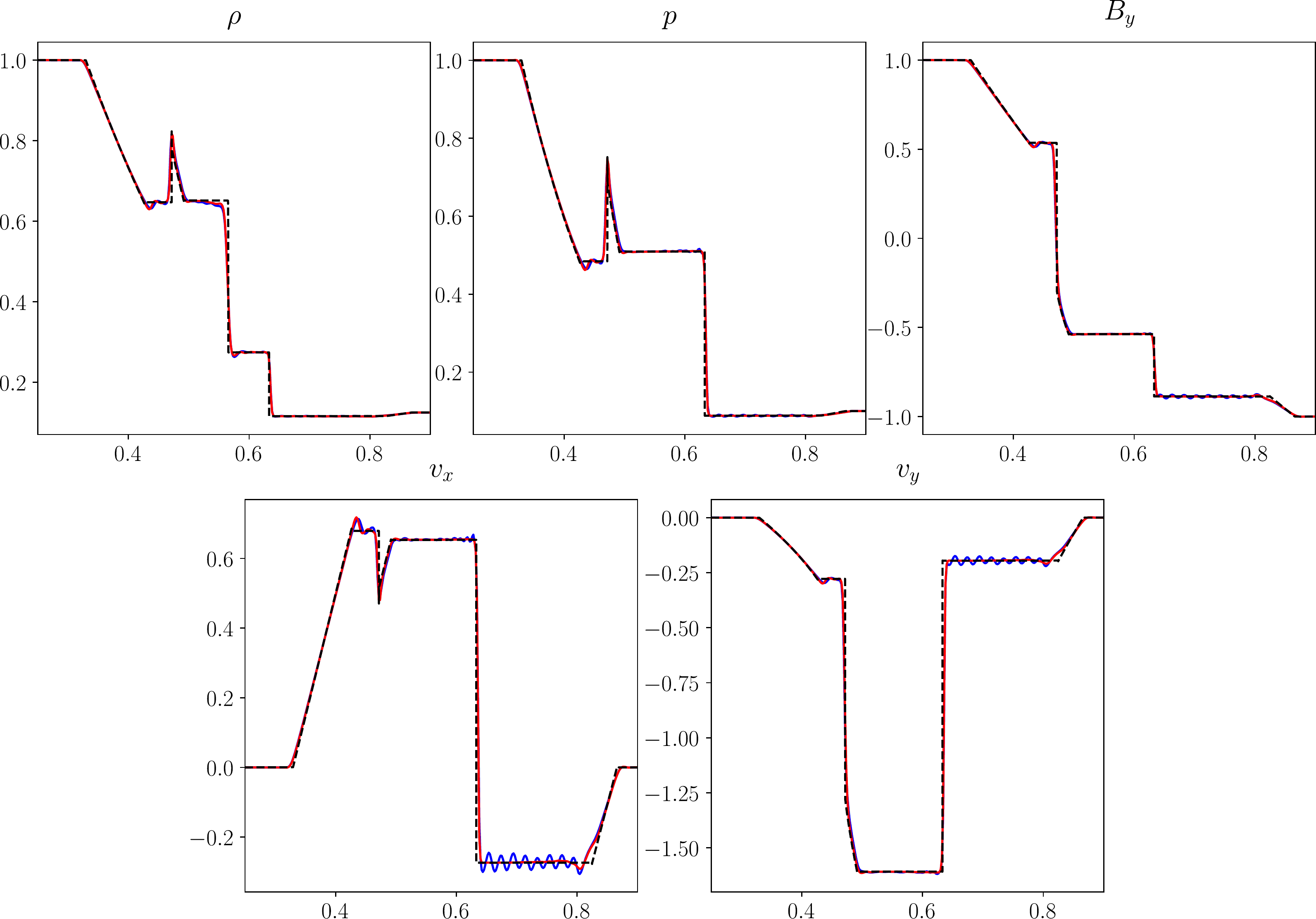}
\caption{ \label{fig:briowu} Solution of the Brio-Wu Riemann problem at $t=0.1$, resolution $512$, CWENO4 scheme. Red and blue: using respectively global and individual smoothness indicators. Black, dashed: with global smoothness indicators at resolution $65536$. Five waves separated by constant states can be seen, from the left to the right: a fast rarefaction wave and a slow compound wave both moving to the left, a contact discontinuity, a slow shock and a fast rarefaction wave, all three moving to the right. }
\end{center}
\end{figure}

One-dimensional Riemann problems serve as standard benchmarks in computational MHD and allow to evaluate the robustness of the numerical method as well as how good it can resolve discontinuities, shocks and appropriate types of waves. In reference \cite{RUJ95}, a lot of Riemann problems are proposed for testing several physical aspects. We present here the results obtained for the classical Brio-Wu shock tube test \cite{BRW88,RUJ95}. The initial conditions for $(\rho,v_x,v_y,v_z,p,B_x,B_y,B_z)$ consist of the constant states $(1,0,0,0,1,0.75,1,0)$ on the left ($x \leq 0.5$) and $(0.125,0,0,0,0.1,0.75,-1,0)$ on the right ($x>0.5$) of the domain $x \in [0,1]$.

In this test, open boundary conditions are used.
Figure \ref{fig:briowu} shows the solution at time $t=0.1$ and provides a comparison between the use of individual/global smoothness indicators. As mentioned at the end of section \ref{sec:WENO}, the use of global smoothness indicators reduces strongly the amount of oscillations.

\subsubsection{Orszag-Tang vortex}
\label{sec:OT}

\begin{figure}
\centering
\includegraphics[width=\textwidth]{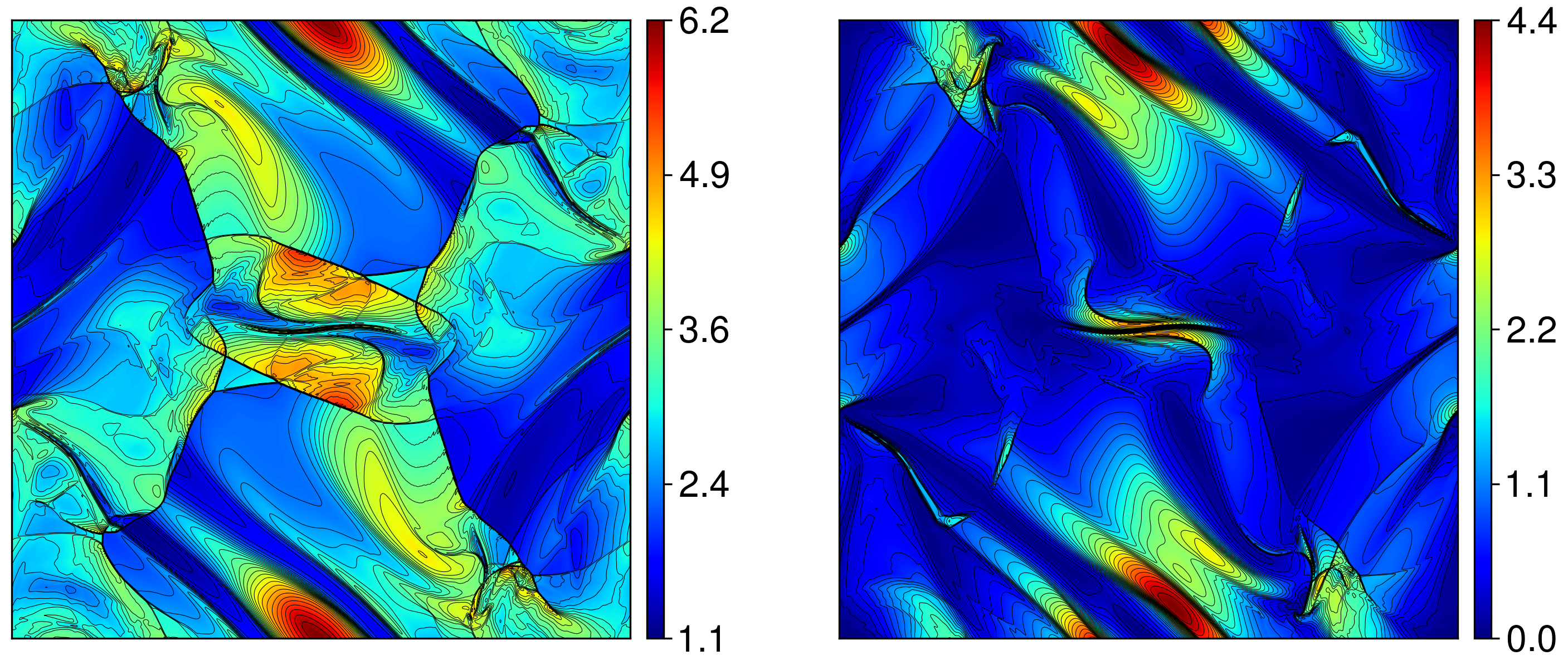}

\begin{tabular}{ cccccccc }
 \hline
\multicolumn{2}{c}{resolution} & $32^2$ & $96^2$ & $288^2$ & $544^2$ & $800^2$ & $1056^2$ \\
\hline
\multirow{2}{*}{$t=0.1$} & $\delta U_{mean}$ & 1.56 $10^{-2}$ & 3.96 $10^{-4}$ & 3.84 $10^{-6}$ & 2.55 $10^{-7}$ & 4.96 $10^{-8}$& 1.54 $10^{-8}$ \\
                        & EOC      & -               & 3.34             & 4.22 & 4.27 & 4.24 & 4.21 \\
\hline
\multirow{2}{*}{$t=0.5$} & $\delta U_{mean}$ & 2.43 $10^{-1}$ & 8.79 $10^{-2}$ & 2.77 $10^{-2}$ & 1.36 $10^{-2}$ & 8.44 $10^{-3}$& 5.75 $10^{-3}$ \\
                        & EOC      & -               & 0.91             & 1.04 & 1.11 & 1.24 & 1.39 \\
\hline
\end{tabular}

\caption{ \label{fig:OTV} Orszag-Tang vortex plots at instant $t=0.5$, resolution $1056^2$, CWENO4 scheme (top), left: density, right: magnetic pressure. Bottom: convergence of errors and EOC at different instants for the CWENO4 scheme.}

\end{figure}

The Orszag-Tang vortex was first studied in the framework of incompressible flows \cite{ORT79} and has become a standard test to check the robustness of compressible MHD solvers (\cite{TOT00,ZIE04,SGA08,MTB10,ROM16} among many others). From smooth initial conditions, many shock structures emerge and develop turbulent dynamics. The initial conditions are given in the computational domain $(x,y) \in [0,1]^2$ by $\rho=\gamma^2$, $\vv=\Big(-\sin(2\pi y),\sin(2\pi x),0\Big)$, $\vB=\Big(-\sin(2\pi y),\sin(4\pi x),0\Big)$ and $p=\gamma$.

Figure \ref{fig:OTV} shows a convergence test done at two instants: $t=0.1$ when the solution is still smooth and at $t=0.5$ with fully developed shocks. Since the exact solution is not known, the point values located at the cell-centers of the lowest resolution run ($32^2$) are compared to a high resolution outcome ($3232^2$) considered as the reference. We can see that even though the convergence order is consistent with a fourth-order scheme at t=0.1, it goes to 1 at later times when the convergence errors are dominated by the shocks, which is consistent with Godunov's theorem \cite{GOD59}. 

\subsubsection{Decaying supersonic MHD turbulence}

\begin{figure}
\begin{center}
\includegraphics[width=\textwidth]{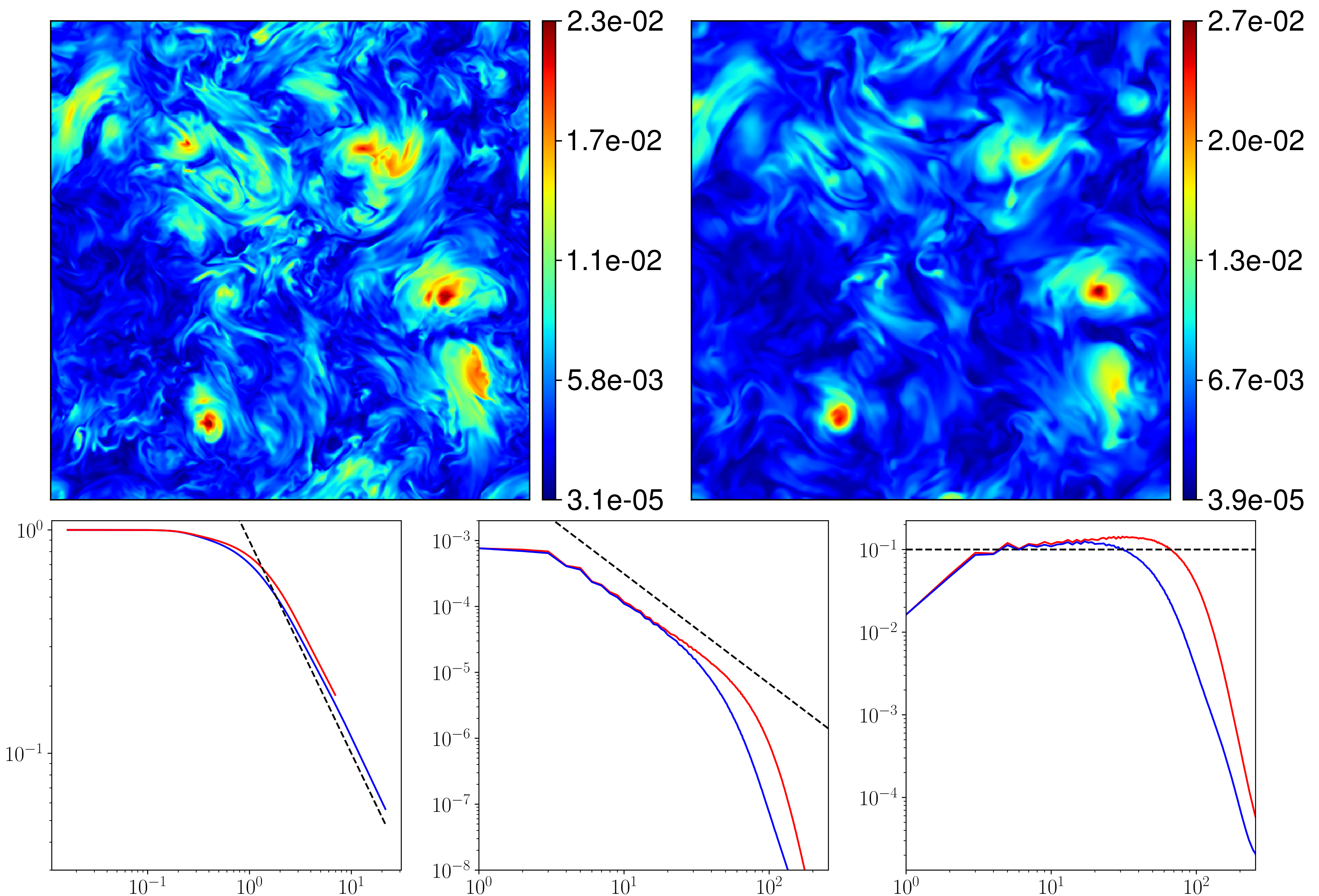}
\caption{ \label{fig:turbdecay} Decaying supersonic $512^3$ MHD turbulence plots comparing the CWENO4FB (red) and the TVD2 (blue) schemes. $(a)$ Slice of kinetic plus magnetic energy at $t \approx 9.84t_E$, CWENO4FB scheme, $(b)$ the same, but with the TVD2 scheme. $(c)$ Log-log plot of the kinetic plus magnetic energy as a function of time (in $t_E$ unit), normalized by its initial value. The dotted line is a guide to the eye and corresponds to the power-law $t^{-0.95}$. $(d)$ Specific kinetic plus magnetic energy spectra at $t \approx 9.84t_E$, the dotted line is a guide to the eye with power-law $k^{-5/3}$. $(e)$ Same spectra but compensated with $k^{5/3}$.}
\end{center}
\end{figure}

As an example of application where the fallback approach is needed, we present here a 3D decaying turbulence experiment at resolution $512^3$. In order to avoid the heating of the plasma when kinetic and magnetic energy are dissipated, we use here the isothermal equation of state (see section \ref{sec:equations} for the implications on the solver) with $c_s=0.1$ on a box $(x,y,z)\in [0,1]^3$. The initial velocity and magnetic field components are generated in Fourier space by giving each mode with wavenumber $\vk$ an amplitude proportional to $\exp(-k^2/(2k_0^2))$ with $k_0=4$, similarly to \cite{BIM00}. The fields are generated in such a way that kinetic helicity, magnetic helicity and cross helicity are initially 0 and the initial velocity field is solenoidal. The fields are then transformed into configuration space and normalized so as to give a Root Mean Squared Mach number of $M=2$ and an initial ratio of magnetic to kinetic energy of 1. The density is initially $\rho=1$ everywhere. 

Without a-priori fallback approach, the CWENO4 scheme gives a negative density shortly after about 0.5 large eddy turnover time. Using the TVD2 scheme as the lower-order method and the thresholds $\tau^{HO}=1, \tau^{LO}=2$ allow the simulation to run without giving negative densities at any time, while the number of reconstructions using a flattener $w^f$ strictly lower than 1 (see section \ref{sec:fallback}) peaked at only about $6/1000$ around the first large eddy turnover time. The large eddy turnover time $t_E$ was computed by considering the instant when the dissipation is maximal. We denote this scheme by CWENO4FB in the following and compare it with TVD2.

At the very beginning of the decay (less than $0.2t_E$), some magnetic energy is rapidly transformed into kinetic energy, accelerating the fluid. The velocity field also causes energy to be stored in a potential form, as density fluctuations, which account for about 10-20$\%$ of the kinetic energy during the run. At $t=t_E$, the full length of scales received energy from larger eddies down to the dissipation scales: the energy dissipation is then maximal and a state of fully-developped turbulence starts to emerge. The energy of the system decays then, following a power-law with an exponent close to 1 (Fig. \ref{fig:turbdecay}c), similarly to the observations made in the nonhelical incompressible case \cite{BIM00}. During the run, the maximum local Mach number peaked at around $7.4$, shortly after $t=t_E$, the maximum local density peaked at about 51.6, the minimum at around 0.01. 

Figure \ref{fig:turbdecay} show comparisons between the CWENO4FB and TVD2 schemes: we can see that the CWENO4FB scheme presents clearly small scale structures at $t\approx 9.84 t_E$ that are smeared out by the second-order scheme. This is also visible when considering the specific kinetic plus magnetic energy Fourier spectra: a Kolmogorov-like inertial range with an exponent close to $-5/3$ is visible for both schemes, but obviously longer for the CWENO4FB scheme. However, when compensating these spectra with $k^{5/3}$, the exponent of the CWENO4FB scheme looks slightly flatter than the one of the TVD2 scheme. 
This well-known effect is probably due to the hyperdiffusive characteristics of the numerical scheme and
can be alleviated by adding explicit physical diffusion terms for viscosity and resistivity.
\section{Final remarks}
In the recent decade significant
progress has been made with regard to higher-order finite-volume algorithms.
This contribution outlines the necessary efforts to go beyond the standard and widespread implementation of
second-order simulation algorithms for ideal magnetohydrodynamics in the conservative finite-volume framework. The gains in precision, detail and numerical consistency with ideal Euler 
dynamics are substantial. The numerical cost remains well acceptable and even decreases dramatically as 
higher-order can offer the same quality as lower-order ones at strongly reduced numerical resolution. 
As the elimination of
significant amounts of numerical dissipation reveals new difficulties that have been covered-up by numerical
diffusion before, strategies of regularizing the solution in the presence of strong shocks are proposed.
Due to the strong reduction of numerical dissipation, schemes that do not strictly preserve the magnetic-field
solenoidality but rely on advection and diffusion might not be as efficient in higher-order simulations. Fortunately, the constrained transport approach
represents a viable solution whose implementation is described here, as well.
In summary, at the current state of numerical developments it seems to be highly appropriate and rewarding
to consider moving existing simulation codes to higher-order schemes.

The contributions of P.-S. Verma and O. Henze to the present work are cordially acknowledged. JMT acknowledges support
by the Berlin International Graduate School in Model and Simulation
based Research (BIMoS).

	\bibliographystyle{abbrv} 
	\bibliography{biblio}

\end{document}